\documentclass[10pt]{article}
\usepackage[english]{babel}
\usepackage[T1]{fontenc}
\usepackage{latexsym}
\usepackage{amssymb}
\usepackage{stmaryrd}
\usepackage{amsmath}
\usepackage{paralist} 
\usepackage[dvips]{graphicx,epsfig}
\usepackage[thmmarks, amsmath]{ntheorem}
\usepackage{color}

\theoremheaderfont{\bf}
\theoremseparator{.}

\theoremseparator{.}
\newtheorem{thm}{Theorem}

\theoremseparator{.}

\theoremseparator{.}

\theoremseparator{.}
\newtheorem{lem}{Lemma}

\theoremseparator{.}
\theorembodyfont{\upshape}
\newtheorem{rem}{Remark}

\newcommand{\halmos}{\vspace{3mm} \hfill \mbox{$\Box$}\\[2mm]}
\newcommand{\conv}[1]{\textrm{conv}}
\newcommand{\var}[1]{\textrm{Var}}

\begin{document}
\begin{center}
{\Large \textbf{On convex hull of Gaussian samples}}\\
\vspace{0.5cm}
{\large Yu. Davydov\footnote{University of Lille 1, France}
}
\end{center}

\vspace{1.5cm}

\noindent\textbf{Abstract:} {\small

Let $X_i = \{X_i(t), \;t \in T\}$ be i.i.d. copies of a centered Gaussian process\\ $X = \{X(t),\;\; t \in T\,\}$ with values in $\mathbb{R}^d$ defined on  a separable metric space $T.$\\ It 
is supposed that $X$ is bounded.
We consider the asymptotic behaviour of convex hulls
$$
W_n = \conv \{\{\,X_1(t),\ldots, X_n(t),\;\;t \in T\}
$$
and show that with probability 1
$$
\lim_{n\rightarrow \infty} \frac{1}{\sqrt{2\ln n}}\,W_n = W
$$
(in the sense of Hausdorff distance), where the limit shape $W$ is defined by the covariance structure of $X$:
$\;\;W = \conv {}\{K_t, \;t\in T\}, \;\; K_t $ being the concentration ellipsoid of $X(t).$

The asymptotic behavior of the mathematical expectations $Ef(W_n),$ where $f$ is an homogeneous functional is also studied
.}\\

\noindent\emph{Key-words:} Gaussian process, Gaussian sample, convex hull, limit theorem.

\section{Introduction}
Let $T$ be a separable metric space. Let $X_i = \{X_i(t), \;t \in T\}$ be i.i.d. copies of a centered Gaussian process $X = \{X(t),\; t \in T\}$ with values in $\mathbb{R}^d.$ Assume that $X$ has a.s. bounded paths and consider the convex hulls
\begin{equation}
\label{W_n}
W_n = \conv \{\{\,X_1(t),\ldots, X_n(t),\;\;t \in T\}.
\end{equation}
We are studying the existence of a limit shape for the sequence $\{W_n\}.$

Our work is motivated by recent papers  \cite{RFMC, MCRF} inspired by an interesting implication in ecological context in estimating the home range of a herd of animals with population size $n.$
Mathematical results of these articles consist in exact computation of a mean perimeter $L_n$ and area $A_n$ of $W_n$ in the case when $d=2$ and $X$ is a standard Brownian motion on $T= [0,1].$
It was shown that
\begin{equation}
\label{perim_area}
L_n \sim 2\pi\sqrt{2\ln n},\;\;\;\; A_n \sim 2\pi\ln n,\;\;\; n\rightarrow \infty.
\end{equation}
The relation between $L_n$ and $A_n$ being the same as the relation between the perimeter and area of a circle of the radius $\sqrt{2\ln n},$ it
seems credible to suppose that $W_n$ rounds up with the growth of $n.$ Our aim is to show that this phenomenon really occurs for all bounded
Gaussian processes. Our main result (Theorem 1) establishes the existence with probability 1 of the limit
\begin{equation}
\label{asympt}
\lim_{n\rightarrow \infty} \frac{1}{\sqrt{2\ln n}}\,W_n = W
\end{equation}
(in the sense of Hausdorff distance) and gives the complete description of the limit set $W$ which is natural to call  {\it limit shape} for
convex hulls $W_n.$ In particular case of standard Brownian motion on $[0,1]$ the set $W$ coincides with the unit ball $B_d(0,1)$ of $\mathbb{R}^d.$

An interesting consequence of (\ref{asympt}) is that the rate of the growth of the convex hulls $W_n$ is the same for all bounded Gaussian processes.

The proof for continuous Gaussian processes may be easily deduced from the known results concerning the asymptotic of Gaussian samples
(see \cite{AK, G}), but in general case one needs
an independent demonstration.

Let us remark in addition that if $T$ is a singleton, $T=\{t_0\},$ and $d=1$, then the process $X$ is simply a real random variable and $W_n$ is the segment\\ $[\max \{\,X_1,\ldots, X_n,\,\}\,,\, 
\min \{\,X_1,\ldots, X_n,\,\}].$ It means that in some sense our study is closely connected with the classical theory of extrema.

\section{Asymptotic behavior of $W_n$}
\subsection{Notation}

$B_d(0,1),\;\;S_d(0,1)$ are respectively unit ball and unit sphere of $\mathbb{R}^d.$

\noindent
$\langle \cdot, \cdot\rangle$ is the scalar product in $\mathbb{R}^d.$

\noindent
$\mathcal{K}(B)$ is the space of compact convex subsets of a Banach space $\mathbb B$ provided with Hausdorff distance $\rho_{\mathbb B}$ :
$$
\rho_{\mathbb B}(A,B) = \max\{\inf\{\,\epsilon \;|\;A \subset B^\epsilon\},\;\; \inf\{\,\epsilon \;|\;B \;\subset A^\epsilon\}\},
$$
$A^\epsilon$ is the open $\epsilon$-neighbourhood of $A$.

\noindent
We set $\mathcal{K}^d = \mathcal{K}(\mathbb{R}^d)$ and $\rho = \rho_{{\mathbb R}^d}.$

\noindent
$\mathcal{M}_A(\theta),\;\;\theta \in S_d(0,1),$ is a support function of a set $A\in \mathcal{K}^d:$
$$
\mathcal{M}_A(\theta)= \sup_{x\in A}\langle x, \theta\rangle, \;\;\;\theta \in S_d(0,1).
$$
$T$ is a separable metric space.

\noindent
$\mathbb{C}(T)$ is the space of continuous functions on $T$ with uniform norm.

\noindent
$X = \{X(t),\; t \in T\}$ is a separable bounded centered Gaussian process with values in $\mathbb{R}^d.$

\noindent
$R_t$ is the covariance matrix of $X(t).$

\noindent
$K_t$ is the ellipsoid of concentration of $X(t):$
$$
K_t = \{x\in \mathbb{R}^d\;|\;\langle R_t^{-1}x, x\rangle \leq 1\}.
$$

\noindent
Finally we set
\begin{equation}
\label{lim_shape}
W = \conv \{\{K_t, \;t\in T\}.
\end{equation}

\subsection{Limit shape}

\begin{thm} \label{as_conv}
${}$

1) Let $X = \{X(t),\; t \in T\}$ be a bounded centered Gaussian process with values in $\mathbb{R}^d.$ Let $(X_i)$ be a sequence of i.i.d. copies of $X$ and $W_n$ be
the convex hull defined by (\ref{W_n}).

 Then with probability 1
\begin{equation}
\label{thm1-1}
 \frac{1}{\sqrt{2\ln n}}\,W_n \;\;
 \stackrel{\mathcal{K}^d}{\longrightarrow}\;\; W.
\end{equation}

2) If $T$ is compact and $X$ is continuous, then a.s.
\begin{equation}
\label{thm1-2}
\rho\left(\frac{1}{\sqrt{2\ln n}}\,W_n,\;\;W\right) = o\left(\frac{1}{\sqrt{\ln n}}\right).
\end{equation}
\end{thm}

\begin{rem} \label{th1}
It is not difficult to see that the support function  $\mathcal{M}_W$ of the limit shape $W$ admits the following representation
$$
\mathcal{M}_W(\theta) = \sigma(\theta),
$$
where
$$
\sigma^2(\theta) = \sup_{t\in T} \langle R_t\theta,\, \theta\rangle,\;\;\; \theta \in S_d(0,1).
$$
\end{rem}
The examples below show that in concrete cases the identification of $W$ is not very complicated.

\begin{rem} \label{th1-2}
For non-centered processes the relation (\ref{thm1-1}) remains the same whe\-reas (\ref{thm1-2}) must be replaced by 
\begin{equation}
\label{thm1-3}
\rho\left(\frac{1}{\sqrt{2\ln n}}\,W_n,\;\;W\right) = O\left(\frac{1}{\sqrt{\ln n}}\right).
\end{equation}
\end{rem}

\subsection{Asymptotic behavior of moments}
Let $f: \mathcal{K}^d \rightarrow \mathbb{R}^1$ be a continuous
positive increasing homogeneous function of degree $p$ , that is
\vspace{5pt}




$f(A) \geq 0 \;\;\; \forall A \in \mathcal{K}^d;$
\vspace{5pt}

$f(A_1) \leq f(A_2) \;\;\;\forall A_1\subset A_2,\;\, A_1,A_2 \in \mathcal{K}^d;$
\vspace{5pt}

$f(cA) = c^pf(A),\;\; \forall \;c\geq 0,\; \forall A \in \mathcal{K}^d.$ \vspace{5pt}







\begin{thm} \label{esper}
Let $f$ be a function with the properties described above. Then, under hypothesis of Theorem 1
\begin{equation}
\label{thm2}
Ef\left(\frac{1}{\sqrt{2\ln n}}\,W_n\right) \rightarrow f\left(W\right).
\end{equation}
\end{thm}

\begin{rem} \label{th2}
${}$

1) This theorem gives in particular the asymptotic behavior for mean values of all reasonable geometrical characteristics of $W_n$ (such as volume, surface measure, diameter,\ ...).

2) By replacing $f$ with $f^m, m>0,$ we get the  asymptotic behavior of higher order moments
\begin{equation}
\label{remthm2}
Ef^m\left(\frac{1}{\sqrt{2\ln n}}\,W_n\right) \rightarrow f^m\left(W\right).
\end{equation}
\end{rem}

\subsection{Examples}

{\bf Brownian motion.} Let $X$ be a standard $d$-dimensional Brownian motion on $T=[0,1].$
Then $K_t = \sqrt{t}B_d(0,1)$ and the limit shape is $W= B_d(0,1).$
In particular, Theorem 2 gives for $d=2$  the relations (\ref{perim_area}).
\vspace{5pt}

\noindent
{\bf Self-similar processes.} Let $X= \{X(t),\;t\in \mathbb{R}_+\}$ be a Gaussian centered self-similar process (SSP) with values in $\mathbb{R}^d.$ It means that for some $\alpha>0$ the processes
$$
\{X(at), \;t\in \mathbb{R}_+\},\;\;\;\{a^{\alpha}X(t), \;t\in \mathbb{R}_+\}
$$
have the same law for any $a>0.$

If we suppose that $X$ is bounded, then we can apply our Theorem 1 to the restriction of $X$ on $[0,1].$ As $X(t) \stackrel{\mathcal{D}}{=} t^{\alpha}X(1),$ we have
$K_t = t^\alpha K_1$ which gives $W=K_1.$

This conclusion is available in  particular when $X$ has the stationary increments: indeed, in this case the process is continuous as for any $\theta \in S_d(0,1)$ the process $\langle X(t),\, \theta\rangle$ is a fractional Brownian motion (FBM). 
\vspace{5pt}

\noindent
{\bf Fractional Brownian Bridge.} Now let us suppose that the coordinates of the process
$Y(t)=\{Y_1(t),\ldots,Y_d(t)\}$ are independent FBM's:
$$
EY_i(t)=0,\;\;\;r(t,s):= EY_i(t)Y_i(s)=\frac{1}{2}(t^{2H} +s^{2H} -|t-s|^{2H}).
$$
 Then the conditional
process related to the condition $Y(1)=0$, which can be called {\it Fractional
  Brownian Bridge,} coincides in distribution with the process
$$
X(t)=\{X_1(t),\ldots,X_d(t)\},\;\;\; X_i(t)=Y_i(t)-r(t,1)Y_i(1),\;\;t\in [0,1].
$$
It is clear that $K_t=\sigma(t)B_d(0,1),$ where $\sigma^2(t) =
t^{2H} -\frac{1}{4}(t^{2H}+1-|1-t|^{2H})^2.$ The function $\sigma^2$ reaches its
maximum at $t=\frac{1}{2}$ and $\sigma^2_{\mathrm {max}}=
\frac{1}{2^{2H}}-\frac{1}{4}.$ Finally we see that $W=\sigma_{\mathrm {max}}B_d(0,1).$

\section{Proofs}

{\bf Proving  Theorem 1.} The theorem is a consequence of two following lemmas.

\begin{lem}
\label{lem1}
Let $Y$ be a r.v. such that for all $\gamma <\frac{1}{2}$
$$
E\exp{\{\gamma Y^2\}} < \infty.
$$
Let $(Y_k)$ be a sequence of independent copies of $Y.$ Then with probability 1
$$
\limsup_n \frac{1}{\sqrt{2\ln n}}\max\{Y_1,\ldots,Y_n\} \leq 1.
$$
\end{lem}

\begin{lem}
\label{lem2}
Let $Y = \sup_T X(t),$ where $X = (X(t), \;t\in T)$ is a centered bounded Gaussian process with
$\sup_T \var {}X(t) = 1$ and let $(Y_n)$ be a sequence of independent copies of $Y.$
 Then with probability 1
$$
\liminf_n \frac{1}{\sqrt{2\ln n}}\max\{Y_1,\ldots,Y_n\} \geq 1.
$$
\end{lem}
\vspace{5pt}

\noindent
 {\usefont{T1}{cmr}{m}{sc}
 \selectfont
  Proof of  Theorem 1, first part.}
Fix $\theta \in S_d(0,1).$ Define $\pi_\theta : \mathbb{R}^d\rightarrow \mathbb{R}^1$ by $\pi_\theta(x) = \langle \theta, x\rangle, \;\;x \in \mathbb{R}^d,$ and set
$$
\pi_\theta(W_n) = [m_n^{(\theta)}, M_n^{(\theta)}],
$$
where
$$
m_n^{(\theta)}= \min_{i\leq n}\{m_{i,\theta}\}, \;\; \;\;
m_{i,\theta}= \inf_{t\in T}\langle \theta, X_i(t)\rangle,
$$
$$
M_n^{(\theta)}=\max_{i\leq n}\{M_{i,\theta}\},\;\; \;\;
M_{i,\theta}= \sup_{t\in T}\langle \theta, X_i(t)\rangle.
$$
Since the paths of $X$ are bounded, due to the well known result of Fernique - Marcus and Shepp \cite{F, MSh} we have that
$$
E\exp{\{\gamma {M_n^{(\theta)}}^2\}} < \infty
$$
for all $\gamma < \frac{1}{2\sigma^2(\theta)},$ where $\sigma^2(\theta)= \sup_{t\in T} \var{} \langle \theta, X_i(t)\rangle.$
\vspace{5pt}

Then by Lemma 1 with probability 1
\begin{equation}
\label{appl1lem1}
\limsup_n \frac{1}{\sqrt{2\ln n}} M_n^{(\theta)}\leq \sigma(\theta).
\end{equation}
On the other hand, by Lemma 2
\begin{equation}
\label{appl1lem2}
\liminf_n \frac{1}{\sqrt{2\ln n}} M_n^{(\theta)}\geq \sigma(\theta), \;\;\;{\textrm a.s.}
\end{equation}
Therefore
\begin{equation}
\label{lem1-2max}
\lim_n \frac{1}{\sqrt{2\ln n}} M_n^{(\theta)} = \sigma(\theta), \;\;\;{\textrm a.s.}
\end{equation}
By the same arguments
\begin{equation}
\label{lem1-2min}
\lim_n \frac{1}{\sqrt{2\ln n}} m_n^{(\theta)} = -\sigma(\theta), \;\;\;{\textrm a.s.}
\end{equation}
It means that for any $\theta \in S_d(0,1)$ with probability 1
\begin{equation}
\label{conv_proj}
\pi_\theta\left( \frac{1}{\sqrt{2\ln n}} W_n\right)\; \longrightarrow\;\;
[- \sigma(\theta), \;\;\sigma(\theta)].
\end{equation}
Let $\{e_i,\; i= 1,\ldots,d\}$ be a basis of $\mathbb{R}^d.$ Consider the parallelepiped $C_n$ defined by the orthogonal projections of $\frac{1}{\sqrt{2\ln n}} W_n$ onto coordinate axes.
The relation (\ref{conv_proj}) implies that
$$
C_n \longrightarrow \prod_{i=1}^d [- \sigma(e_i), \;\;\sigma(e_i)], \;\;\;{\textrm a.s.}
$$
Since $\frac{1}{\sqrt{2\ln n}} W_n \subset C_n,$ the sequence $\{\frac{1}{\sqrt{2\ln n}} W_n\}$ is bounded, and hence relatively compact, in $\mathcal{K}^d.$

Due to the natural isometry between $(\mathcal{K}^d,\,\rho)$ and $\mathbb{C}(S_d(0,1))$ it follows from this that the sequence $\{\mathcal{M}_n(\theta),\;\; \theta \in S_d(0,1)\}$ of support functions of the sets $\{\frac{1}{\sqrt{2\ln n}} W_n\}$ is a.s. relatively compact in the space $\mathbb{C}(S_d(0,1)).$ Let $\Theta$ be a countable dense subset of $S_d(0,1).$ Using the relation
(\ref{lem1-2max}) we see that with probability 1 for all $\theta \in \Theta$
$$
\mathcal{M}_n(\theta) \rightarrow \sigma(\theta).
$$
Together with relative compactness this shows that almost surely the sequence $\{\mathcal{M}_n(\cdot)\}$ has a unic limit point. Then the same is true for $\{\frac{1}{\sqrt{2\ln n}} W_n\},$ and Remark 1 concludes the proof of the first part.
\vspace{8pt}

\noindent
 {\usefont{T1}{cmr}{m}{sc}
 \selectfont
  Proof of  Theorem 1, second part.} Now we can consider the processes $X$ and $X_k$ as random elements of the separable Banach space $\mathbb{B}=\mathbb{C}(T).$ By Theorem 2.1. of \cite{G} with probability 1
\begin{equation}
\label{goodman}
\rho_{\mathbb{B}}\left(\widetilde{W}_n,\;\;\sqrt{2\ln n}\,\widetilde{W}\right) = o(1),
\end{equation}
where $\widetilde{W}_n = \conv{}_{\mathbb{B}}\{X_1,\ldots,X_n\}$ and $\widetilde{W}$ is the ellipsoid of concentration  of $X.$ Let $\varphi : \mathbb{B} \rightarrow \mathcal{K}^d$ be defined by
$$
\varphi(x) = \conv {}\{\,x(t), \;t\in T\}.
$$
It is clear that $\varphi(\widetilde{W}_n)= W_n,\;\;\;\varphi(\widetilde{W})= W,$ and it is easy to check that the map $\varphi$ is Lipschitzian:
$$
\rho(\varphi(x),\;\varphi(y)) \leq \rho_{\mathbb{B}}(x,y),\;\; x,y \in \mathbb{B}.
$$
Therefore (\ref{thm1-2}) follows directly from (\ref{goodman}).
$\halmos$

{\bf Proofs of Lemmas 1-2.}
\vspace{5pt}

\noindent
 {\usefont{T1}{cmr}{m}{sc}
 \selectfont
  Proof of  Lemma 1.} Let $s>0.$ Setting
$$
Z_n =  \frac{1}{\sqrt{2\ln n}}\max\{Y_1,\ldots,Y_n \}
$$
  our assumption implies
  $$
  P\left\{Z_n \geq 1+s \right\}\;\; \leq  \;\;
  nP\{Y\geq (1+s)\sqrt{2\ln n}\}\leq
$$
$$
\hspace{-40pt}\leq
 \;\;\frac{n E\exp{\{\gamma Y^2\}}}{\exp{\{\gamma (1+s)^22\ln n\}}}\;\; =\;\; C(\gamma)n^{-\delta},
  $$
where $\delta = 1 - 2\gamma(1+s)^2 >0$ if $\gamma > \frac{1}{2(1+s)^2}.$

This inequality shows that the series
$$
\sum_m P\left\{Z_{m^a} \geq 1+s \right\}
$$
is summable if $a > 1/\delta.$
We apply the Borel--Cantelli lemma and since $s$ is arbitrary, we find that
$$
\limsup_m Z_{m^a} \leq 1.
$$
As for $k \in [m^a, (m+1)^a)$
$$
Z_k \leq Z_{(m+1)^a}\sqrt{\frac{\ln (m+1)}{\ln m}},
$$
we get
$$
\limsup_n Z_{n} \leq 1.
$$
\vspace{-28pt}

$\halmos$

\noindent
 {\usefont{T1}{cmr}{m}{sc}
 \selectfont
  Proof of  Lemma 2.} Let $0<s<1$ and let $t_0 \in T$ be chosen so that $\sigma_0^2 :=\var {}X(t_0) > s.$ We use the same notation $Z_n$ for $\frac{1}{\sqrt{2\ln n}}\max\{Y_1,\ldots,Y_n \}.$  We have
$$
P\left\{Z_n \leq s \right\}\;\;= P\left\{Y \leq s \sqrt{2\ln n}\right\}^n\;\; = \;\;
 F(s\sqrt{2\ln n})^n\;\;\;
 $$
 $$
 \leq  \;\;\exp{\{-n(1-F(s\sqrt{2\ln n}))\}},
$$
where $F$ is the distribution function of $Y.$

Note  that
$$
1-F(x) = P\{\sup_T X(t) > x\}\;\; \geq \;\;P\{X(t_0) > x\}\;\; \geq\;\;Cx^{-1}\exp{\left\{-\frac{x^2}{2\sigma_0^2}\right\}},
$$
which shows that
$$
P\left\{Z_n \leq s \right\}\;\;\leq \exp{\left\{-C(\ln n)^{-\frac{1}{2}}n^{1-\frac{s^2}{\sigma_0^2}}\right\}}.
$$
It means that the series $\sum_n P\left\{Z_n \leq s \right\}$ is summable, and by
applying the Borel--Cantelli lemma we finish the proof. $\halmos$
\vspace{10pt}

{\bf Proving Theorem 2.} The proof is based on the following lemma completing the information given by Lemma 1.

\begin{lem}
\label{lem3}
Let $Y$ be a r.v. such that for some $\gamma > 0$
$$
E\exp{\{\,\gamma Y^2\}} < \infty.
$$
Let $(Y_i)$ be a sequence of independent copies of $Y.$
Then for any $k\in {\mathbb N}$
$$
\sup_n E \left(\frac{1}{\sqrt{2\ln n}}\max\{\,Y_1,\ldots,Y_n\}\right)^k < \infty.
$$
\end{lem}
\noindent
 {\usefont{T1}{cmr}{m}{sc}
 \selectfont
  Proof of  Theorem 2.}
Due to the continuity of $f$ and the convergence (\ref{thm1-1}) the result will follow from the uniform integrability of the family $\left\{f\left(\frac{W_n}{\sqrt{2\ln n}}\right)\right\}.$
Using the notation from the proof of the first part of Theorem 1, we set
$$
m_n = \min_{i=1,\ldots,d} m_n^{e_i},\;\;\;M_n = \max_{i=1,\ldots,d} M_n^{e_i},\;\;\;
D_n = \max\{- m_n,M_n\},
$$
$$
L_n = [-D_n, D_n]^d.
$$
Since $W_n \subset L_n,$ we have
$$
f\left(\frac{W_n}{\sqrt{2\ln n}}\right) \;\;\leq \;\;f\left(\frac{L_n}{\sqrt{2\ln n}}\right)
\;\;=\;\; \left(\frac{D_n}{\sqrt{2\ln n}}\right)^p f([-1,1]^d),
$$
hence it is sufficient to state that for all $p>0$
\begin{equation}
\label{supE}
\sup_n E\left(\frac{D_n}{\sqrt{2\ln n}}\right)^p \;\;< \infty.
\end{equation}
 The latter relation  follows directly from Lemma 1 , and the
theorem is proved. $\halmos$

\noindent
 {\usefont{T1}{cmr}{m}{sc}
 \selectfont
  Proof of  Lemma 3.} Let
$$
Z_n =  \frac{1}{\sqrt{2\ln n}}\max\{\,Y_1,\ldots,Y_n \}.
$$
Denote by $F,\; F_n$ the distribution functions of $Y$ and $Z_n$ respectively. By Markov inequality and by the assumption of Lemma
\begin{equation}
\label{lem3_1}
1-F_n(x) \leq nP\{Y\geq x\sqrt{2\ln n}\} \;\;\leq  An^{1-2\gamma x^2}.
\end{equation}
Hence for $a=\frac{1}{\sqrt{2\gamma}}$
\begin{equation}
\label{lem3_2}
 E(Z_n)^k = \int_0^\infty x^{k-1}(1-F_n(x))dx \;\;\leq
 a^k\;+ \;kA\int_a^\infty x^{k-1}n^{1-2\gamma x^2}dx.
\end{equation}
As for $x\geq \frac{1}{\sqrt{2\gamma}}$
$$
n^{1-2\gamma x^2} \leq \exp{\{-\gamma x^2\ln n\}},
$$
 we find that
$$
\limsup_n \int_a^\infty x^{k-1}n^{1-2\gamma x^2}dx \;\;=\;\; 0,
$$
and we get from (\ref{lem3_2})
$$
\limsup_n E(Z_n)^k\;\;\leq \;\;a^k,
$$
which completes the proof. $\halmos$

\section{Concluding remarks}

\hspace{11pt} 1. It is clear that the result of the second part of Theorem 1 is still available if the space
$\mathbb{R}^d$ is replaced by a separable Banach space.

On the contrary, the similar question about  the first part of Theorem 1 and  about Theorem 2 is more delicate: their proofs are essentially based on the
 compactness of bounded subsets of $\mathbb{R}^d$ and it is not clear how to handle this obstacle in the infinite-dimensional case.

\vspace{10pt}

2. The second interesting question is about the character of approach of $\frac{W_n}{\sqrt{2\ln n}}$ to $W$. Is it true that $b_n \rho\left(\frac{W_n}{\sqrt{2\ln n}}\,,\;W\right)$ does converge in law to some limit for some choice of normalizing constants $b_n$?

The same question can be also asked for  the processes \\ $\{b_n(M_n(\theta) - M(\theta))\,,\;\;\theta \in S^{d-1}\},$
where $M_n(\theta),\; M(\theta)$ are respectively the support functions of $\frac{W_n}{\sqrt{2\ln n}}$ and $W.$ It seems that the recent paper \cite{K} may be useful in this context.
\vspace{10pt}

3. What  can we say on the behaviour of $W_n$ in non-Gaussian
case? It is more or less clear that the convergence a.s. must be
replaced by the weak one and the  normalizing constants will be
transformed from logarithmic to power ones. Indeed, if $X$ is a
vector in $\mathbb{R}^m$ with regulary varying distribution and $(X_i)$ is a
sequence of i.i.d. copies of $X$, then it is  well known that the
point processes
$$
\beta_n = \sum_{i=1}^n \delta_{\left\{\frac{X_i}{n^{1/\alpha}}\right\}}
$$
converge weekly to some Poisson point process $\Pi_\alpha.$

It follows immediately from this that
\begin{equation}
\label{stable}
\frac{W_n}{n^{1/\alpha}} \Longrightarrow  \conv{}(\Pi_\alpha).
\end{equation}
This fact is still available for a much more general case when $X$ is a random element of an abstract convex cone $\mathbb K$ (see \cite{DMZ}) provided $X$ satisfies the condition of regular variation (condition (4.5) in \cite{DMZ}). 
Hence the convergence (\ref{stable}) may  be considered as a ``regular varying'' analog of the second part of  Theorem 1.
The main difficulty now is how to check the condition of regular variation for concrete situations.

\vspace{10pt}

{\bf Acknowledgments.} The author wishes to thank M. Lifshits 
for his interest to this work and useful discussions and V. Paulauskas for stimulating remarks,
 as well as all participants of working seminar on Stochastic Geometry of the university Lille 1 for their support.



\begin{thebibliography}{00}

\vspace{10pt}





\bibitem{AK} A. De Acosta and J. Kuelbs, Limit theorems for moving averages,\\ Z. Wahrsch. verw. Gebiete, 64 (1983), pp. 67--123.

\bibitem{DMZ} Yu. Davydov, I. Molchanov and S. Zuyev,
Strictly stable distributions on convex cones, EJP,
 13 (2008),  11, pp. 259--321.

\bibitem{F}  X. Fernique, Régularité de processus gaussiens, Invent. Math., 12 (1971), pp. 304--320.

\bibitem{RFMC} J. Randon-Furling, Satya N. Majumdar and A. Comptet, Perimeter and Area of the Convex Hull of $N$ Planar Brownian Motions, preprint, ArXiv:0907.0921v1, 6 Jul 2009.

\bibitem{G} V. Goodman, Characteristics of normal samples,
Ann. Probab. 16 (1988), 3, pp. 1281--1290.

\bibitem{K} Z. Kabluchko, M. Schlauter and L. de Haan, Stationary max-stable fields associated to negative definite functions, Ann. Probab. 37 (2009), 5, pp. 2042--2065.

\bibitem{L} M. Lifshits, Gaussian Random Functions, Kluwer (1995), 337p.

\bibitem{MCRF} Satya N. Majumdar, A. Comptet and J. Randon-Furling, Random convex hulls and extreme value statistics, preprint, ArXiv:0912.0631v1, 3 Dec 2009.

\bibitem{MSh} M. B. Marcus and L. A. Shepp, Sample behavior of Gaussian processes, Proc. Sixth Berkeley Symp. Math. Statist. Prob., 2 (1971), pp. 423--442.







\end{thebibliography}
\end{document}